\documentclass[a4paper,12pt]{article}
\usepackage[a4paper,top=4.0cm,bottom=2.0cm,left=2.0cm,right=2.0cm]{geometry}

\usepackage{natbib}
\bibpunct{(}{)}{;}{d}{,}{,}

\usepackage{amsthm,amsmath,amsfonts,mathrsfs}
\usepackage{epsfig,graphicx}

\title{A mollified ensemble Kalman filter\thanks{Universit\"at Potsdam, 
Institut f\"ur Mathematik, Am Neuen Palais 10, D-14469 Potsdam, Germany}}
\author{
Kay Bergemann
\and
Sebastian Reich}

\begin{document}

\maketitle

\begin{abstract}
It is well recognized that discontinuous analysis increments of sequential data assimilation
systems, such as ensemble Kalman filters, might lead to spurious high frequency adjustment
processes in the model dynamics. Various methods have been devised 
to continuously spread out the analysis increments over a fixed time interval centered about analysis 
time. Among these techniques are nudging and incremental analysis updates (IAU). 
Here we propose another alternative, which may be viewed as a hybrid of nudging and IAU and which 
arises naturally from a recently proposed continuous formulation of the ensemble Kalman analysis 
step. A new slow-fast extension of the popular Lorenz-96 model is introduced to demonstrate the
properties of the proposed mollified ensemble Kalman filter.

\end{abstract}

\section{Introduction}

Given a model in form of an ordinary differential equation 
and measurements at discrete instances in time, data assimilation attempts to find the best possible
approximation to the true dynamics of the physical system under consideration
\citep{sr:evensen}. Data assimilation 
by sequential filtering techniques achieves such an approximation by discontinuous-in-time 
adjustment of the model dynamics due to available measurements. While optimality of
the induced continuous-discrete filtering process can be shown for linear systems, 
discontinuities can lead to unphysical readjustment processes under the model dynamics
for nonlinear systems and under imperfect knowledge of the error probability distributions. See, e.g.,
\cite{sr:houtekamer05,sr:kepert09} in the context of operational models as well as 
\cite{sr:neef06,sr:oke07} for a study of this phenomena under simple model problems.
These observations have led to the consideration of data assimilation systems that seek to
incorporate data in a more ``continuous'' manner. In this paper we focus on a novel continuous
data assimilation procedure based on ensemble Kalman filters. Contrary to widely used incremental analysis 
updates (IAU) of \cite{sr:bloom96}, which first do a complete analysis step to then distribute
the increments evenly over a given time window, our approach performs the analysis step itself over 
a fixed window. Our novel filter technique is based on the continuous formulation of the Kalman 
analysis step in terms of ensemble members
\citep{sr:br10} and mollification of the Dirac delta function by smooth approximations 
\citep{sr:friedrichs}. The proposed mollified ensemble Kalman (MEnK) filter is described in
Section \ref{sec:MEnKF}. The MEnK filter may be viewed as a ``sophisticated'' form of
nudging \citep{sr:hoke76,sr:macpherson91} with the nudging coefficients being 
obtained from a Kalman analysis perspective instead of heuristic tuning. However, for certain
multi-scale systems, nudging with prescribed nudging coefficients might still be advantageous 
\citep{sr:kalnay09a}. We also point to the closely related work by \cite{sr:stauffer09}, which
proposes a nudging-type implementation of the ensemble Kalman filter 
with perturbed observations \citep{sr:burgers98}. 

To demonstrate the properties of the new MEnK filter, we propose a slow-fast extension
of the popular Lorenz-96 model \citep{sr:lorenz96} in Section \ref{sec:slowfast}. 
Contrary to other multi-scale extensions of the Lorenz-96 model, the fast dynamics is entirely 
conservative in our model and encodes a dynamic balance relation similar to geostrophic 
balance in primitive equation models. Our model is designed to show the generation of 
unbalanced fast oscillations through standard sequential ensemble Kalman filter 
implementations under imperfect knowledge of the error probability distributions. Imperfect knowledge
of the error probability distribution can arise for various reasons. In this paper, we address in particular
the aspect of small ensemble size and covariance localization \citep{sr:houtekamer01,sr:hamill01}.
It is also demonstrated that both IAU as well as the newly proposed MEnK 
filter maintain balance under assimilation with small ensembles and covariance localization. 
However, IAU develops an instability in our model system over longer assimilation cycles, which 
requires a relatively large amount of artificial damping in the fast dynamics. 
We also note that the MEnK filter is cheaper to implement than
IAU since no complete assimilation cycle and repeated model integration need to be performed.
It appears that the MEnK filter could provide a useful alternative to the widely employed
combination of data assimilation and subsequent ensemble re-initialization. See, 
for example, \cite{sr:houtekamer05}.

\section{Mollified ensemble Kalman filter} \label{sec:MEnKF}

We consider models given in form of ordinary differential equations
\begin{equation} \label{ode}
\dot{\bf x} = f({\bf x},t)
\end{equation}
with state variable ${\bf x} \in \mathbb{R}^n$. Initial conditions 
at time $t_0$ are not precisely known and 
we assume instead that
\begin{equation} \label{IC}
{\bf x}(t_0) \sim {\rm N}({\bf x}_0,{\bf B}),
\end{equation}
where ${\rm N}({\bf x}_0,{\bf B})$ denotes an $n$-dimensional 
Gaussian distribution with mean
${\bf x}_0 \in \mathbb{R}^n$ and covariance matrix ${\bf B} \in 
\mathbb{R}^{n \times n}$. We also
assume that we obtain measurements ${\bf y}(t_j) \in \mathbb{R}^{p}$ 
at discrete times $t_j\ge t_0$, $j=0,1,\ldots,M$, subject to 
measurement errors, which are also Gaussian 
distributed with zero mean and covariance matrix ${\bf R} 
\in \mathbb{R}^{p\times p}$, i.e.
\begin{equation} \label{measurement}
{\bf y}(t_j) - {\bf H}{\bf x}(t_j) \sim {\rm N}({\bf 0},{\bf R}).
\end{equation}
Here ${\bf H} \in \mathbb{R}^{p\times n}$ is the (linear) measurement 
operator.

Ensemble Kalman filters \citep{sr:evensen} rely on 
the simultaneous propagation of $m$ independent solutions ${\bf x}_i(t)$, $i=1,\ldots,m$, 
of (\ref{ode}), which we collect in a matrix ${\bf X}(t) \in
\mathbb{R}^{n\times m}$. We can extract an empirical mean
\begin{equation}
\overline{\bf x}(t) = \frac{1}{m} \sum_{i=1}^m {\bf x}_i(t)
\end{equation}
and an empirical covariance matrix
\begin{equation}
{\bf P}(t) = \frac{1}{m-1} \sum_{i=1}^m \left({\bf x}_i(t)-\overline{\bf x}(t)\right)
 \left({\bf x}_i(t)-\overline{\bf x}(t)\right)^T
\end{equation}
from the ensemble.
In typical applications from meteorology, the ensemble size $m$ is much smaller than the
dimension $n$ of the model phase space. We now state a complete ensemble Kalman filter 
formulation for sequences of observations at time instances $t_j$, $j=1,\ldots,M$, and 
intermediate propagation of the ensemble under the dynamics 
(\ref{ode}). Specifically, the continuous formulation
of the ensemble Kalman filter step by \cite{sr:br10} allows for 
the following concise formulation in terms of a single differential equation 
\begin{equation} \label{odeEnKF}
\dot{\bf x}_i =
f({\bf x}_i,t) - \sum_{j=1}^M
\delta(t-t_j)\,{\bf P} \nabla_{{\bf x}_i} {\cal V}_j({\bf X})
\end{equation}
in each ensemble member, where $\delta(\cdot)$ denotes the standard Dirac 
delta function, ${\cal V}_j$ is the potential
\begin{equation} \label{potential}
{\cal V}_j({\bf X}) = \frac{m}{2} \left\{ S_j(\overline{\bf x}) 
+ \frac{1}{m} \sum_{i=1}^m S_j({\bf x}_i)
\right\} 
\end{equation}
with observational cost function
\begin{equation} \label{S}
S_j({\bf x}) = \frac{1}{2} \left({\bf H}{\bf x} -{\bf y}(t_j)\right)^T
{\bf R}^{-1} \left({\bf H}{\bf x} -{\bf y}(t_j)\right).
\end{equation}
One may view (\ref{odeEnKF}) as the original ODE (\ref{ode}) 
driven by a sequence of impulse like contributions due to observations. 
See the Appendix for a derivation of (\ref{odeEnKF}).

It should be noted that (\ref{odeEnKF}) is equivalent to a standard Kalman filter
(and hence optimal) for linear models and the number of ensemble members $m$ larger 
than the dimension of phase space $n$. This property does no longer hold under nonlinear
model dynamics and it makes sense to ``mollify'' the discontinuous analysis adjustments
as, for example, achieved by the popular IAU technique \citep{sr:bloom96}. 
In the context of (\ref{odeEnKF}), a most natural mollification is 
achieved by replacing (\ref{odeEnKF}) with
\begin{equation} \label{mollyEnKF}
\dot{\bf x}_i =
f({\bf x}_i,t) -  \sum_{j=1}^M
\delta_\epsilon(t-t_j)\,{\bf P} \nabla_{{\bf x}_i} {\cal V}_j({\bf X}),
\end{equation}
where 
\[
\delta_\epsilon(s) = \frac{1}{\epsilon} \psi(s/\epsilon),
\]
$\psi(s)$ is the standard hat function 
\begin{equation} \label{hat}
\psi(s) = \left\{ \begin{array}{ll} 1-|s| & \mbox{for}\,\,|s|\le 1,\\
0 & \mbox{else} \end{array} \right. 
\end{equation}
and $\epsilon>0$ is an appropriate parameter. 
The hat function (\ref{hat}) could, of course, be replaced by another B-spline. We note that
the term mollification was introduced by \cite{sr:friedrichs} to denote families of 
compactly supported smooth functions $\delta_\epsilon$ which approach the 
Dirac delta function $\delta$ in the limit $\epsilon \to 0$.  Mollification
via convolution turns non-standard functions (distributions) into smooth functions.
Here we relax the smoothness assumption and allow for any non-negative, compactly 
supported family of functions that can be used to approximate the Dirac delta function.
A related mollification approach to numerical time-stepping methods has been proposed
in the context of multi-rate splitting methods for highly oscillatory
ordinary differential equations  in \cite{sr:garcia-archilla97tmi,sr:izaguirre99jcp}.

Formulations (\ref{odeEnKF}) and (\ref{mollyEnKF}) are based on square-root filter 
formulations of the ensemble Kalman filter \citep{sr:tippett03}. Ensemble inflation 
\citep{sr:andand99} is a popular techniques to stabilize the performance of such
filters. We note that ensemble inflation can be interpreted as adding a
$\theta \left( {\bf x}_i - \overline{\bf x} \right)$ term to the right hand side of
(\ref{mollyEnKF}) for an appropriate parameter value $\theta >0$. We note that a related
``destabilizing'' term appears in $H_\infty$ continuous filter formulations for linear systems. 
See, for example, \cite{sr:simon} for details.  The link to $H_\infty$ filtering suggests 
alternative forms of ensemble inflation and a systemtic approach to make ensemble Kalman
filters more robust. Details will be explored in a forthcoming publication. 

Localization \citep{sr:houtekamer01,sr:hamill01} is another popular technique to enhance the
performance of ensemble Kalman filters. Localization can easily be implemented into the mollified
formulation (\ref{mollyEnKF}) and leads to a modified 
covariance matrix $\widetilde{\bf P} = {\bf C} \circ {\bf P}$ in (\ref{mollyEnKF}). 
Here ${\bf C}$ is an appropriate localization
matrix and ${\bf C} \circ {\bf P}$ denotes the Schur product of ${\bf C}$ and ${\bf P}$. See
\cite{sr:br10} for details.

We finally note that both  (\ref{odeEnKF}) and (\ref{mollyEnKF}) can be modified to become 
consistent with an ensemble Kalman filter with perturbed observations \citep{sr:burgers98}. This 
requires to adjust the potentials ${\cal V}_j$ and to add a stochastic forcing term.

\section{Algorithmic summary}

Formulation (\ref{mollyEnKF}) can be solved numerically by any standard ODE 
solver. Similar to nudging, there is no longer a strict separation 
between ensemble propagation and filtering. 
However, based on (\ref{mollyEnKF}), ``nudging'' is performed consistently 
with Kalman filtering in the limit $\epsilon \to 0$. An optimal choice 
of $\epsilon$ will depend on the specific applications. 

Using, for example, the implicit midpoint method as the time-stepping method for the
dynamic part of (\ref{mollyEnKF}) and denoting numerical approximations 
at $t_k = k\Delta t$ by ${\bf x}_i^k$, we obtain the following algorithm:
\begin{equation} \label{algorithm}
{\bf x}_i^{k+1} = {\bf x}_i^k + \Delta t \left\{
f({\bf x}_i^{k+1/2},t_{k+1/2}) -  \sum_{j=1}^{M}\alpha_{j}^k {\bf P}^k \nabla_{{\bf x}_i}
{\cal V}_j({\bf X}^k) \right\}, 
\end{equation}
for $i=1,\ldots,m$. 
Here ${\bf X}^k$ is the collection of ensemble approximations
at $t_k$, ${\bf P}^k$ is the resulting covariance matrix, ${\bf x}_i^{k+1/2} = ({\bf x}_i^{k+1} + 
{\bf x}_i^k)/2$ is the midpoint approximation, and the weights $\alpha_{j}^k \ge 0$ are defined by
\begin{equation}
\alpha_{j}^k = \frac{c}{\epsilon} \psi 
\left(\frac{t_k-j \Delta t_{\rm obs}}{\epsilon}\right)
\end{equation}
with the hat function (\ref{hat}) and normalization constant $c>0$ chosen such that
\begin{equation}
 \Delta t \sum_k \alpha_{j}^k  = 1 .
\end{equation}
Note that the summation in (\ref{algorithm}) needs to be performed only over those (measurement) 
indices $j$ with $\alpha_{j}^k \not= 0$.

For the numerical 
experiments conducted in this paper, we set $\epsilon = \Delta t_{\rm obs}/2$ and assume that
$\Delta t_{\rm obs}$ is a multiple of the time-step size $\Delta t$. Under this setting only
a single measurment ${\bf y}(t_j)$ is assimilated at any given time-step $t_k$.
Another interesting choice is $\epsilon = \Delta t_{\rm obs}$ which implies that $\sum_j \alpha_{j}^k$ is 
independent of the time-step index $k$. In this case, two measurements are simultaneously assimilated 
(with different weights) at any given time-step $t_k$.

The main additional cost of (\ref{algorithm}) over standard implementations of
ensemble Kalman filters is the computation of the covariance matrix
${\bf P}$ in each time-step. However, this cost is relatively minor compared to the
cost involved in propagating the ensemble under the model equations (\ref{ode}) in
the context of meteorological applications. We note that, contrary to IAU and the
recent work by \cite{sr:stauffer09}, no additional model propagation steps are required.

\section{A slow-fast Lorenz-96 model} \label{sec:slowfast}

We start from the standard Lorenz-96 model \citep{sr:lorenz96,sr:lorenz98}
\begin{equation} \label{Lorenz}
\dot{x}_l = (x_{l+1} - x_{l-2})x_{l-1} - x_l + 8
\end{equation}
for $l=1,\ldots,40$ with periodic boundary conditions, i.e.,
$x_{-1} = x_{39}$, $x_0 = x_{40}$, and $x_{41} = x_1$.
We note that the energy 
\begin{equation}
E_{\rm Lorenz} = \frac{1}{2} \sum_{l=1}^{40} x_l^2
\end{equation}
is preserved under the ``advective'' contribution
\begin{equation}
\dot{x}_l = (x_{l+1} - x_{l-2}) x_{l-1} = x_{l-1}x_{l+1} - x_{l-2}x_{l-1}
\end{equation}
to the Lorenz-96 model. 

We now derive a slow-fast extension of (\ref{Lorenz}).
To do so we introduce an additional variable $h$, which satisfies
a discrete wave equation, i.e., 
\begin{equation} \label{Wave}
\varepsilon^2 \ddot{h}_l = -h_l + \alpha^2 \left[h_{l+1} 
-2h_l + h_{l-1}\right] ,
\end{equation}
$l = 1,\ldots,40$, where $\varepsilon>0$ is a small parameter implying that
the time evolution in $h$ is fast compared to that of $x$ in
the Lorenz model (\ref{Lorenz}) and $\alpha >0$ is a parameter determining
the wave dispersion over the computational grid. 
We assume again periodic boundary conditions and find that (\ref{Wave})
conserves the energy
\begin{equation}
E_{\rm wave} = \frac{\varepsilon^2}{2} \sum_{l=1}^{40} \dot{h}_l^2 
+ \frac{1}{2} \sum_{l=1}^{40} \left[ h_l^2 + \alpha^2
\left(h_{l+1}-h_l\right)^2 \right].
\end{equation}

\begin{figure}
\begin{center}
\includegraphics[width=0.5\textwidth]{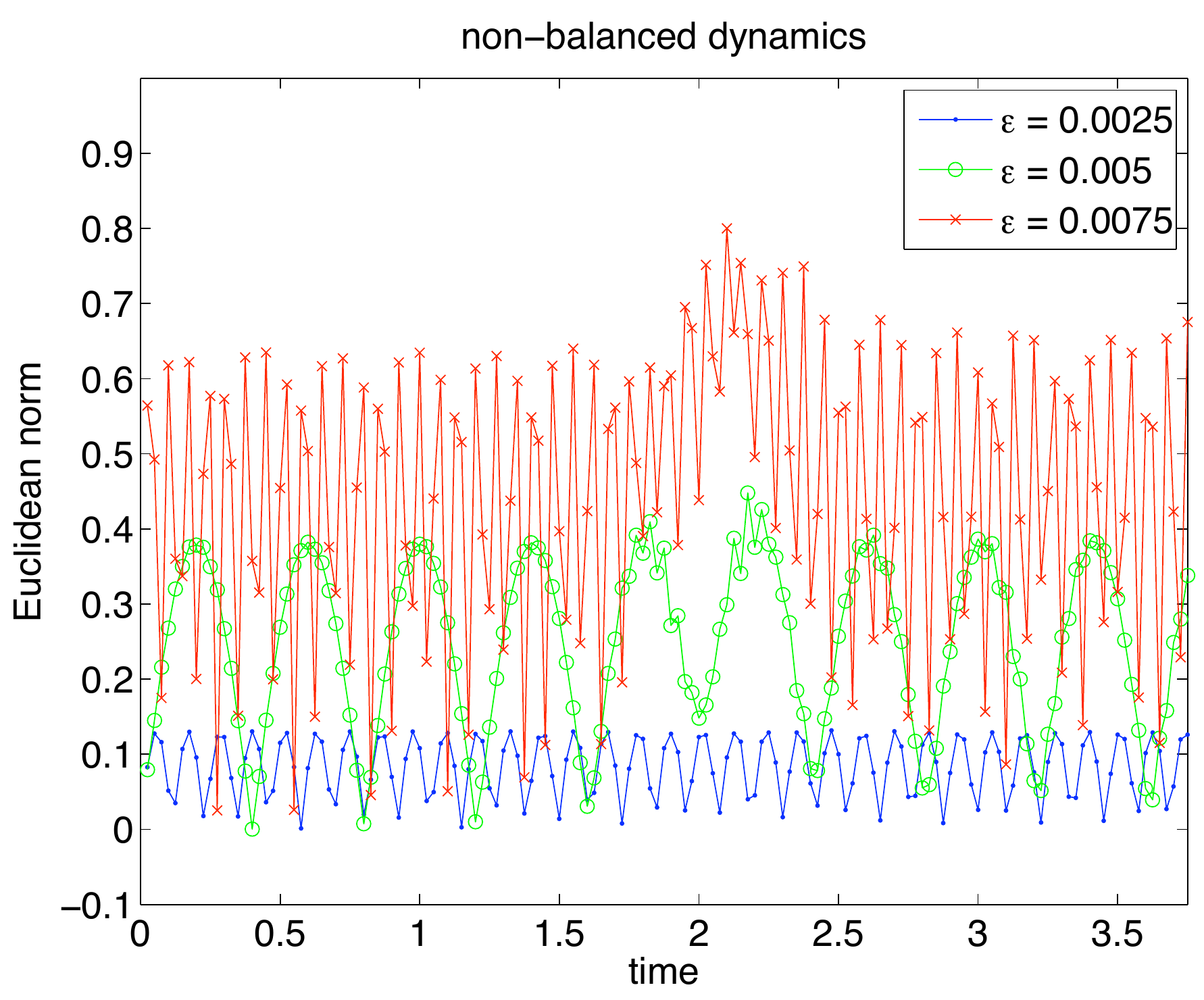}
\end{center}
\caption{Norm of (\ref{norm}) over all grid points 
for several values of $\varepsilon$ as a function of integration 
time. It can be concluded that balance is maintained to high accuracy for sufficiently small values
of $\varepsilon$ except for perturbations induced at initial time. The amplitude of the oscillations 
could be reduced further by more sophisticated initialization techniques (higher order balance 
conditions).}
\label{fig_0}
\end{figure}

\begin{figure}
\begin{center}
\includegraphics[width=0.5\textwidth]{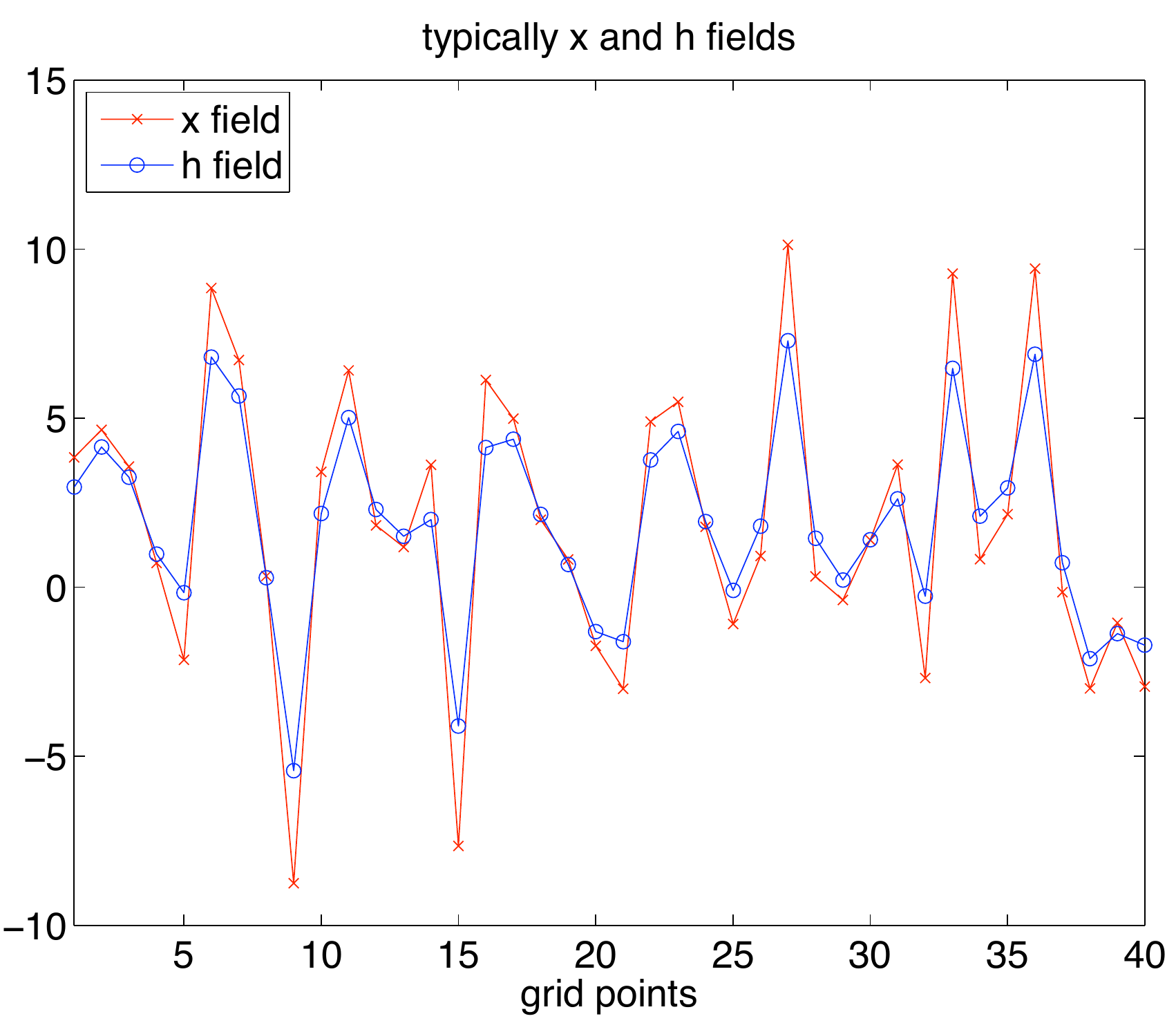}
\end{center}
\caption{Typical balanced $\{x_l\}$ and $\{h_l\}$ fields for the slow-fast
Lorenz model. It can be seen that the balanced $\{h_l\}$ is more regular 
than the $\{x_l\}$ field due to the balance relation (\ref{balance}).}
\label{fig_0b}
\end{figure}

The two models (\ref{Lorenz}) and (\ref{Wave}) are now coupled 
together by introducing an exchange energy term
\begin{equation}
E_{\rm coupling} = - \delta \sum_{l=1}^{40} h_l x_l,
\end{equation}
where $\delta > 0$ characterizes the coupling strength.  
The equations of motion for the
combined system are then defined as
\begin{eqnarray} \label{mixed1}
\dot{x}_l &=& (1-\delta)\left(x_{l+1} - x_{l-2}\right)x_{l-1} 
+ \delta \left( x_{l-1}h_{l+1} - x_{l-2}h_{l-1}\right) - x_l + 8,\\
\varepsilon^2 \ddot{h}_l &=& -h_l + \alpha^2 \left[h_{l+1} 
-2h_l + h_{l-1}\right] +  x_l, \label{mixed2}
\end{eqnarray}
where we have scaled the advective contribution in (\ref{Lorenz})
by a factor of $1-\delta$ to counterbalance the additional contribution
from the coupling term. We note that the pure wave-advection system
\begin{eqnarray} \label{mixed1b}
\dot{x}_l &=&  (1-\delta) \left( x_{l+1} - x_{l-2}\right)x_{l-1}  
+ \delta \left( x_{l-1}h_{l+1} - x_{l-2}h_{l-1}\right),\\
\varepsilon^2 \ddot{h}_l &=& - h_l + \alpha^2\left[h_{l+1} 
-2h_l + h_{l-1}\right] +  x_l \label{mixed2b}
\end{eqnarray}
conserves the total energy
\begin{eqnarray}
H &=& (\delta -1 ) E_{\rm Lorenz} + \delta E_{\rm wave} + E_{\rm coupling}\\
&=& \frac{\delta}{2} \sum_{l=1}^{40} \left\{
\frac{\delta -1}{\delta} x_l^2 + \varepsilon^2 \dot{h}_l^2 +
h_l^2 + \alpha^2\left(h_{l+1}-h_{l}\right)^2 -2x_l h_l \right\}
\end{eqnarray}
independent of the choice of $\delta \in [0,1]$.
We also wish to point to the balance relation
\begin{equation} \label{balance}
x_l = h_l - \alpha^2 \left[h_{l+1} 
-2h_l + h_{l-1}\right]
\end{equation}
in the wave part (\ref{mixed2}). More specifically, if (\ref{balance}) is
satisfied at initial time, the variables $\{h_l\}$ and $\{x_l\}$ will
remain in approximate balance over a fixed time interval
with the deviation from exact balance proportional
to $\varepsilon$ \citep{sr:wirshep00,sr:wir04,sr:cotter06}.
We demonstrate conservation of balance for decreasing values
of $\varepsilon$ in Fig.~\ref{fig_0}, where 
we display the Euclidean norm of
\begin{equation} \label{norm}
{\Delta}_l = x_l - h_l + \alpha^2 \left[h_{l+1} -2h_l + h_{l-1}\right],
\end{equation}
$l=1,\ldots,40$,  as a function of time.
Hence, we may view (\ref{balance}) as defining an approximative
slow manifold on which the dynamics of (\ref{mixed1})-(\ref{mixed2}) 
essentially reduces to the
dynamics of the standard Lorenz-96 model for sufficiently small
$\varepsilon$. We use $\varepsilon = 0.0025$ for our numerical experiments. 
On the other hand, any imbalance introduced through
a data assimilation system will remain in the system since there
is no damping term in (\ref{mixed2}). This is in contrast to other
slow-fast extensions of the Lorenz-96 model, such as 
\cite{sr:lorenz96}, which
also introduce damping and forcing into the fast degrees of freedom. 
We will find in Section \ref{sec_numerics} that some ensemble Kalman filter implementations 
require the addition of an ``artificial'' damping term into (\ref{mixed2}) and
we will use
\begin{equation} \label{mixed2D}
\varepsilon^2 \ddot{h}_l = -h_l + \alpha^2 \left[h_{l+1} 
-2h_l + h_{l-1}\right] +  x_l - \gamma \varepsilon^2 \dot{h}_l
\end{equation}
for those filters with an appropriately chosen parameter $\gamma>0$.

We mention that (\ref{mixed1})-(\ref{mixed2}) may be viewed as an one-dimensional 
``wave-advection'' extension of the Lorenz-96 model under 
a quasi-geostrophic scaling regime with $\varepsilon$ being proportional to the 
Rossby number and ${\rm Bu} = \alpha^2$ being the Burger number \citep{sr:salmon99}. 
In a quasi-geostrophic regime, the Burger number should be of order one and
we set $\alpha = 1/2$ in our experiments. Since a typical time-scale of
the Lorenz-96 model is about $0.025$ units, we derive at an equivalent
Rossby number of $\mbox{Ro} \approx 0.1$ for $\varepsilon = 0.0025$.
Typical $\{ h_l \}$ and $\{ x_l \}$ fields are displayed in Fig.~\ref{fig_0b}.
It should, of course, be kept in mind that 
the Lorenz-96 model is not the discretization of a realistic fluid model. 

The dynamics of our extended model depends on the coupling strength $\delta \in [0,1]$. Since 
the $\{h_l\}$ field is more regular than the $\{x_l\}$ field, stronger coupling leads to a more regular dynamics 
in terms of spatial correlation. We compute the grid-value mean $\overline x = \langle x_l \rangle$ and its variance 
$\sigma = \langle (x_l-\bar x)^2 \rangle^{1/2}$ along a long reference trajectory for $\delta = 0.1$ 
($\bar x = 2.32$, $\sigma = 3.68$), $\delta = 0.5$ ($\bar x = 1.80$, $\sigma = 3.67$), and 
$\delta = 1.0$ ($\bar x = 1.48$, $\sigma = 3.69$), respectively. The variance is nearly identical for 
all three parameter values while the mean decreases for larger values of $\delta$, which implies that the 
dynamics becomes less ``non-linear'' for larger values of $\delta$.

\section{Numerical experiments} \label{sec_numerics}

\begin{figure}
\begin{center}
\includegraphics[width=0.5\textwidth]{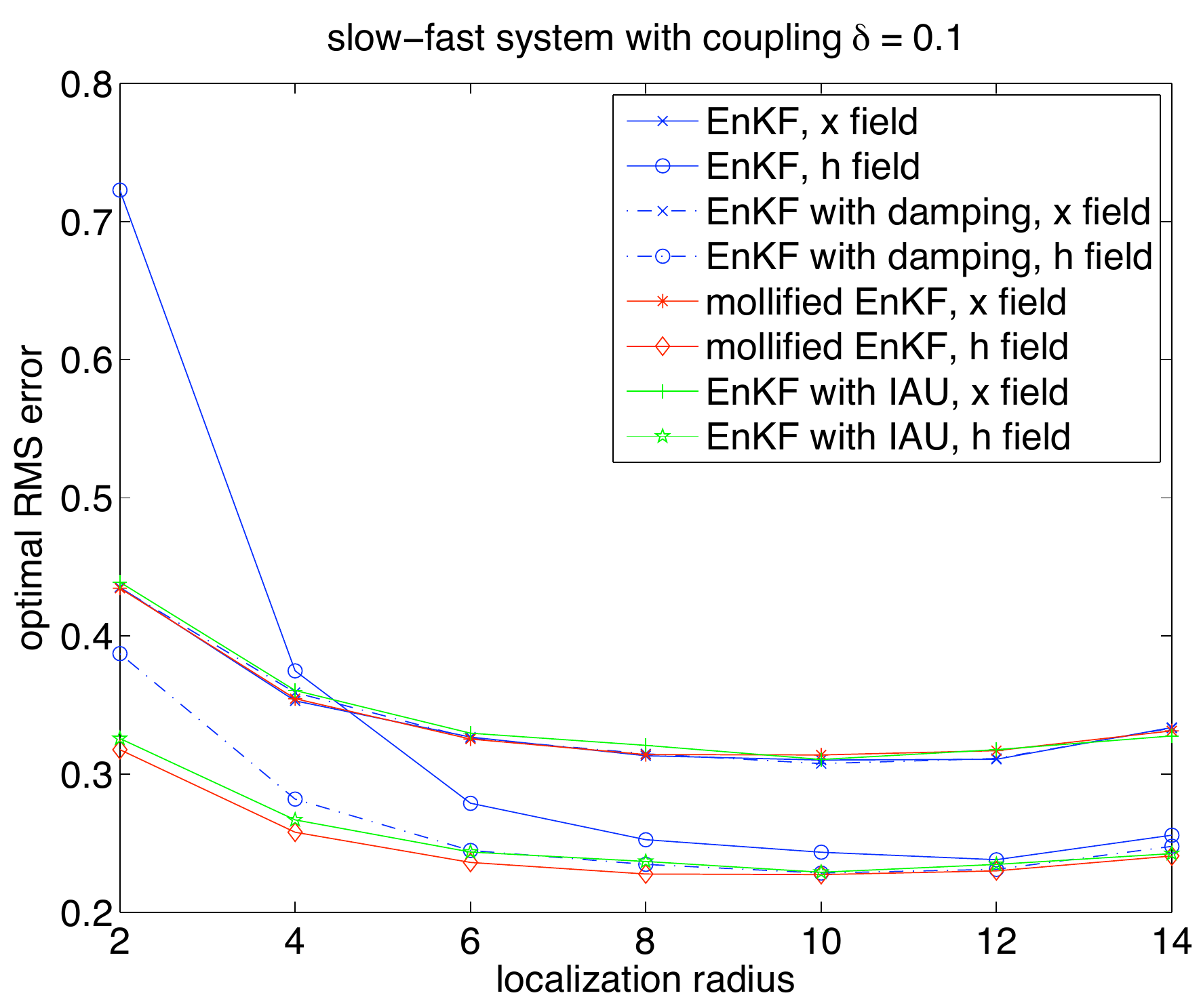}
\end{center}
\caption{Best RMS error for slow-fast Lorenz-96 model with coupling strength
$\delta = 0.1$ using ensembles of size $m = 10$ and $k=20$ observations 
of $\{x_j\}$ taken in intervals of $\Delta t_{\rm obs} = 0.05$ over a total
of 4000 assimilation cycles. We compare a standard ensemble Kalman filter 
(EnKF) implementation to the MEnK filter and a EnKF with IAU. We also 
implement a standard EnKF for the damped
wave equation (\ref{mixed2D}) with damping factor $\gamma = 0.1$.}
\label{fig_1}
\end{figure}

\begin{figure}
\begin{center}
\includegraphics[width=0.5\textwidth]{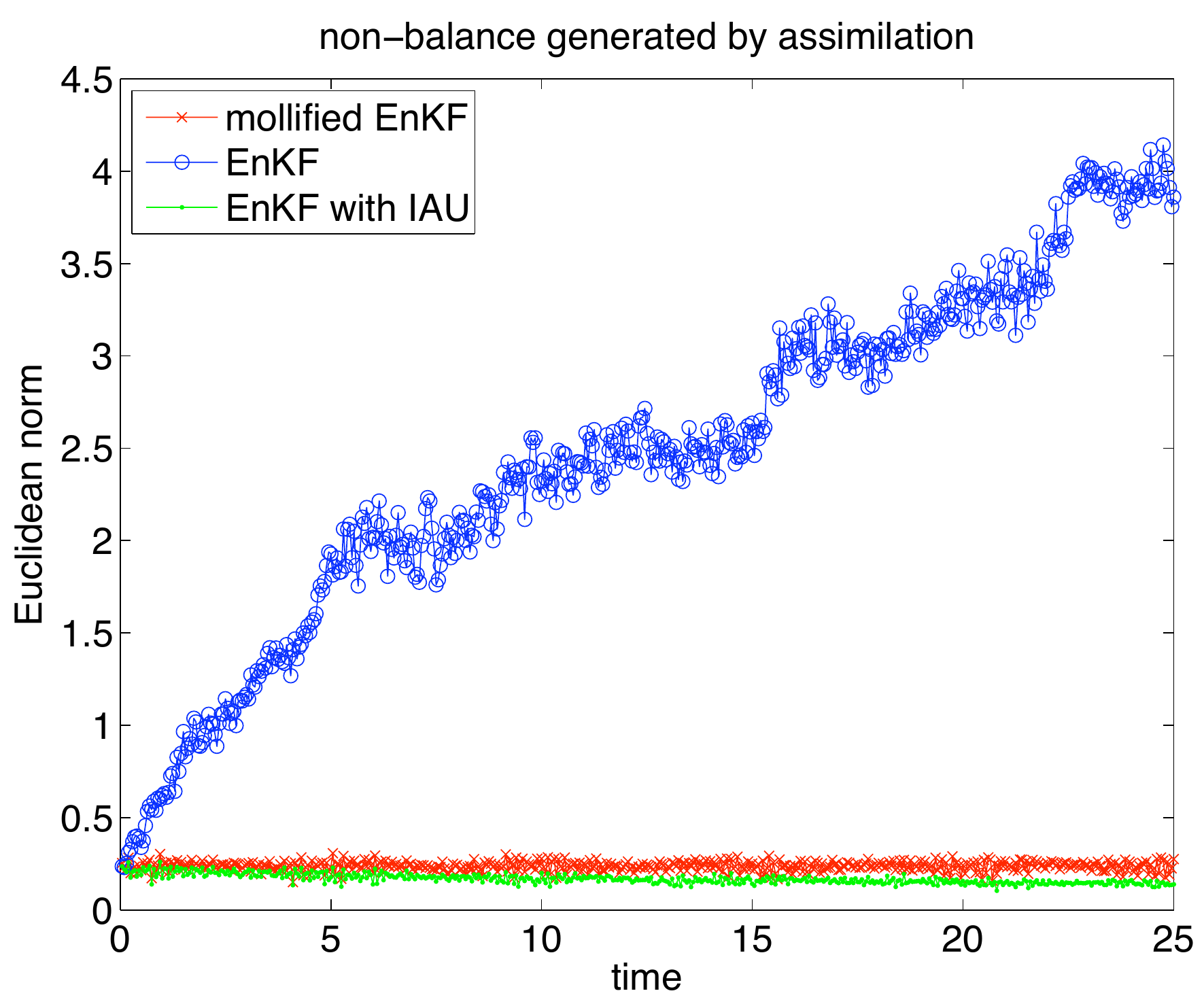}
\end{center}
\caption{Unbalanced dynamics generated through the data assimilation steps
in a standard ensemble Kalman filter (EnKF)
implementation as compared to the relatively low and nearly constant unbalanced wave
activities under the MEnK filter and the IAU EnKF implementation over 500 assimilation steps and 
for coupling constant $\delta = 0.1$ and localization radius $r_0 = 2$.}
\label{fig_3}
\end{figure}

We run the system for a coupling strengths of $\delta = 0.1$ and $\delta = 0.5$, respectively, 
using either a (standard) ensemble Kalman filter implementation based 
on (\ref{odeEnKF}) or an implementation of the mollified formulation 
(\ref{mollyEnKF}). We also recall that the mollified formulation 
(\ref{mollyEnKF}) uses the mollifier (\ref{hat}) with 
$\epsilon = \Delta_{\rm obs}/2$.

We observe $\{x_l\}$ at every second grid point in time-intervals of $\Delta t_{\rm obs} = 0.05$ with
measurement error variance ${\bf R} = {\bf I}_{20}$. The time-step for the 
numerical time-stepping method (a second-order in time, time-symmetric 
splitting method \citep{sr:LeiRei04}) 
is $\Delta t = \Delta t_{\rm obs}/20 = 0.0025$.

All filter implementations use localization \citep{sr:houtekamer01,sr:hamill01} and 
ensemble inflation \citep{sr:andand99}. Localization is performed by multiplying each element 
of the covariance matrix ${\bf P}$ by a distance dependent factor
$\rho_{i,i'}$. This factor is defined by the compactly supported localization function
(4.10) from \cite{sr:gaspari99}, distance 
$r_{l,l'} = \min\{|l-l'|,40-|l-l'|\}$, where $l$ and $l'$ denote the
indices of the associated observation/grid points $x_{l}$ and $x_{l'}$, respectively, and
a fixed localization radius $r_0$. See \cite{sr:br10} for more implementation details. 
Inflation is applied after each time-step to the $\{x_l\}$ components of the ensemble 
members. 

The initial ensemble is generated by adding small random perturbations to a reference 
$\{x_l\}$ field. The associated $\{h_l\}$ ensemble fields are computed from the balance relation
(\ref{balance}). After a spin-up period of 
200 assimilation cycles, we do a total of 4000 assimilation cycles in each 
simulation run and compute the RMS error between the analyzed predictions and the true reference 
solution from which the measurements were generated. RMS errors are computed over a wide range of 
inflation factors and only the optimal result for a given localization radius is reported.

\begin{figure}
\begin{center}
\includegraphics[width=0.5\textwidth]{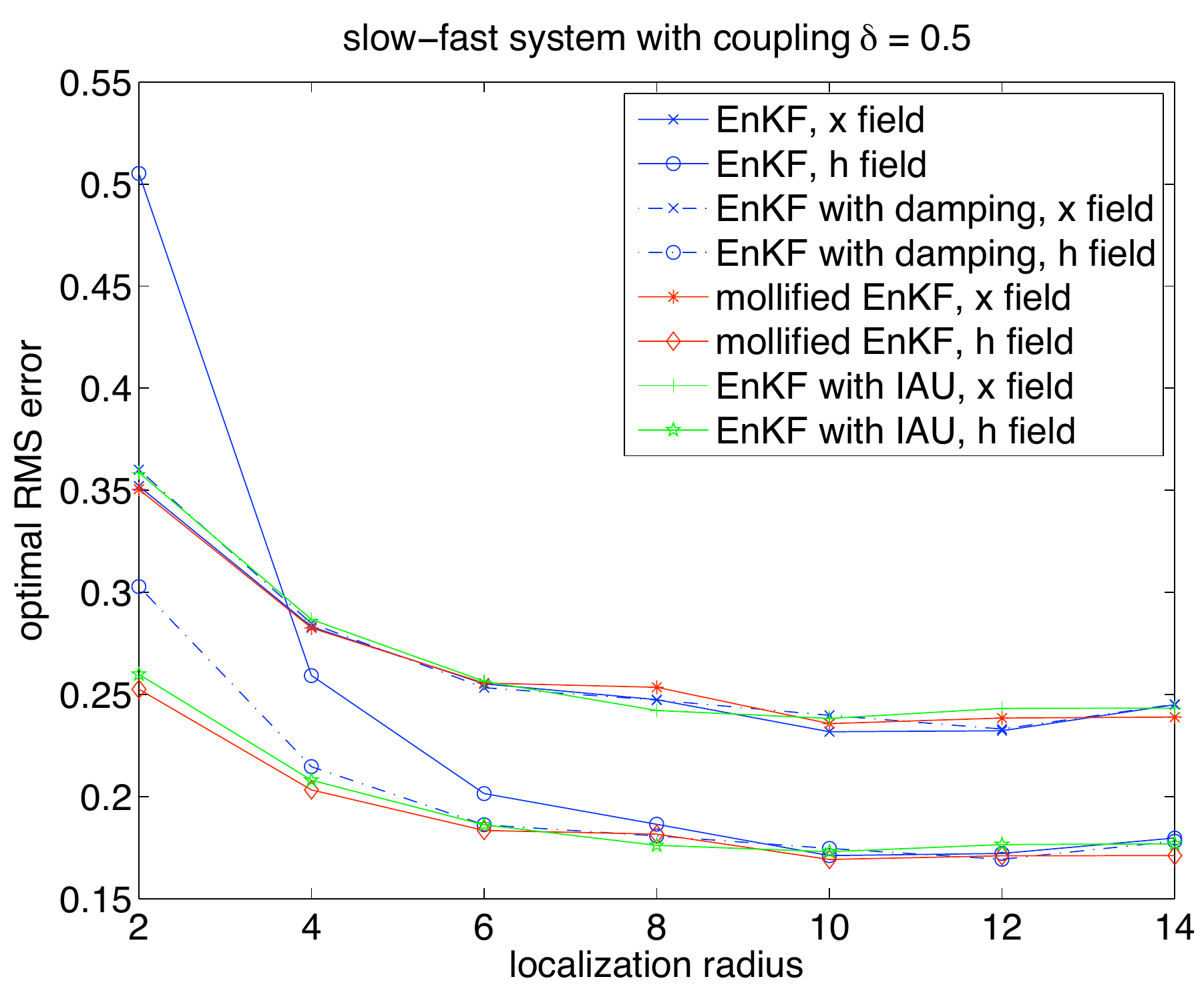}
\end{center}
\caption{Same as Fig.~\ref{fig_1} except that the coupling strength is set to
$\delta = 0.5$.}
\label{fig_2}
\end{figure}

We first compare the performance of a standard ensemble Kalman filter implementation with
the proposed MEnK filter for a coupling strength $\delta = 0.1$. 
It can be seen from Fig.~\ref{fig_1} that the RMS error 
in the $\{h_l\}$ fields is much larger for the standard 
ensemble Kalman filter than for the mollified scheme. The difference is most pronounced
for small localization radii. Indeed, it is found that the standard filter generates 
significant amounts of unbalanced waves. See Fig.~\ref{fig_3}, where we display the Euclidean 
norm of (\ref{norm}) over all ensemble members for the first 500 assimilation steps and 
localization radius $r_0 = 2$. Note that ${\Delta}_l = 0$ at $t=0$ and that no damping is 
applied to the wave part (\ref{mixed2}). We also implement a standard
ensemble Kalman filter for the damped wave equation (\ref{mixed2D}) with $\gamma = 0.1$
which leads to a significant reduction (but not complete elimination) 
in the RMS errors for the $\{ h_l \}$ fields. Qualitatively the same results are obtained for
the stronger coupling $\delta = 0.5$, see Fig.~\ref{fig_2}.

For comparison we implement an IAU along the lines of 
\cite{sr:bloom96,sr:polavarapu04} with the analysis
increments provided by a standard ensemble Kalman filter. Instead of a constant forcing, we use 
the hat function (\ref{hat}) for computing the weights $g(t)$ in \cite{sr:bloom96} to be consistent
with the implementation of the MEnK filter. While it is found that
the resulting IAU ensemble Kalman filter conserves balance well initially (see Fig.~\ref{fig_3}), 
the IAU filter eventually becomes highly inaccurate and/or unstable. We believe
that this failure of IAU over longer analysis cycles (here 4000) is due to an instability
of the algorithm in the fast wave part (\ref{mixed2}). The instability could be caused by
the non-symmetric-in-time nature of the IAU process, i.e., by the fact that one first
computes a forecast, then finds the analysis increments based on the forecast
and available observations, and finally repeats the forecast with assimilation increments included.  
It should also be kept in mind that the IAU of \cite{sr:bloom96}) is based on a standard 3D Var
analysis for given background covariance matrix, which is significantly different from 
an IAU ensemble Kalman filter implementation. A stable implementation of an IAU is achieved by
replacing (\ref{mixed2}) by the damped version (\ref{mixed2D}) with $\gamma = 1.0$. 
Smaller values of $\gamma$ lead to filter divergence. With damping included, the IAU filter behaves
almost identical to the MEnK filter and the standard ensemble Kalman filter with damping
in the wave part (\ref{mixed2D}), see Figs.~\ref{fig_1} and \ref{fig_2}.

\begin{figure}
\begin{center}
\includegraphics[width=0.5\textwidth]{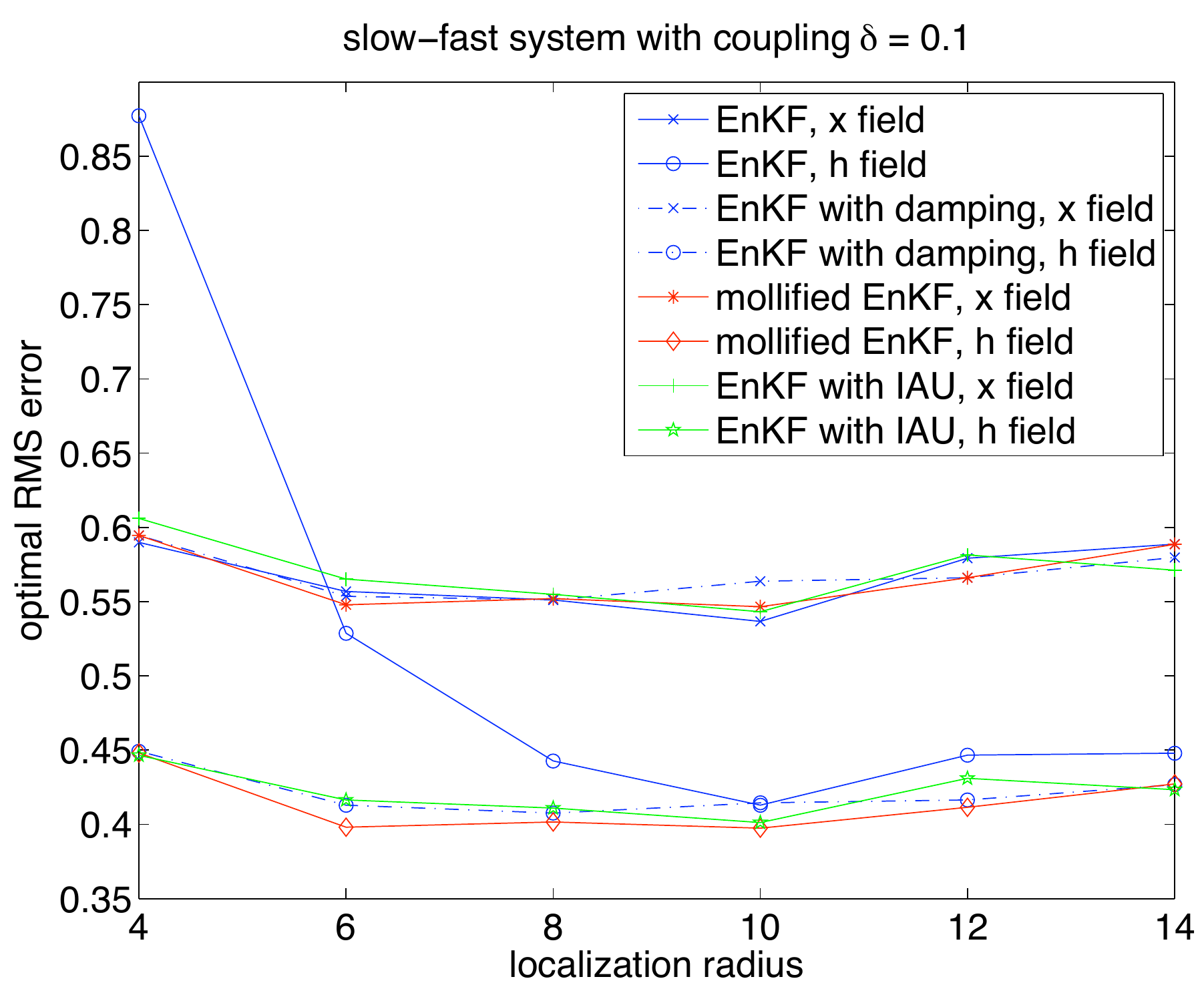}
\end{center}
\caption{Same as Fig.~\ref{fig_1} except that observations of 
$\{(x_j+h_j)/2\}$ are taken in intervals of $\Delta t_{\rm obs} = 0.05$ with 
measurement error variance ${\bf R} = {\bf I}_{20}$ 
over a total of 4000 assimilation cycles.}
\label{fig_2b}
\end{figure}

We also tested the dependence on the type of measurements and found
that observation of the mixed $\{(h_l+x_l)/2\}$ field 
instead of $\{x_l\}$ leads to results which are in qualitative agreement 
with those displayed in Fig.~\ref{fig_1} for all three methods, 
see Fig.~\ref{fig_2b}.

In summary, we find that localization interferes with dynamic balance relations in standard ensemble Kalman 
filter implementations. This result has been well established by a number of previous studies on models of
varying complexity. Here we have demonstrated for a relatively simply model system that strong localization
can be used without the need of additional damping and/or re-initialization provided the assimilation increments 
are distributed evenly over a time-window instead of being assimilated instantaneously. It appears that the proposed
MEnK filter works particularly reliably. It should be kept in mind that the MEnK filter leads to additional costs per
model integration step. But these costs should be relatively small compared to those of the model dynamics itself.

\section{Conclusion}

In this paper, we have addressed a shortcoming of sequential data assimilation systems
which arises from the discontinuous nature of the assimilation increments and the imperfect
statistics under ensemble based filtering techniques. In line with the previously considered
IAU and nudging techniques, we have proposed spreading the analysis through a mollified
Dirac delta function. Mollification has been used  before to avoid numerical
instabilities in multi-rate time-stepping methods for the integration of highly oscillatory
differential equations (see, for example, \cite{sr:garcia-archilla97tmi,sr:izaguirre99jcp}). 
We have demonstrated for a simple slow-fast extension of the Lorenz-96 model 
that the proposed MEnK filter indeed eliminates high-frequency responses in the dynamic model. 
The MEnK filter can be viewed as a nudging technique with the nudging coefficients determined
consistently from ensemble predictions. We note in this context
that a combined ensemble Kalman filter and nudging technique was found beneficial 
by \cite{sr:kalnay09a} for the original slow-fast Lorenz-96 model \citep{sr:lorenz96}.  
The MEnK filter (perhaps also combined with the filtering-nudging approach of 
\cite{sr:kalnay09a}) may provide an alternative to currently used combinations of ensemble 
Kalman filters and subsequent re-initialization. For the current model study it is found that
artificial damping of the unbalanced high-frequency waves stabilizes the performance of a
standard ensemble Kalman filter and leads to results similar to those observed for the
MEnK filter. We also found, in line with previous studies of \cite{sr:houtekamer05,sr:oke07,sr:kepert09}, 
that stronger localization leads to a more pronounced generation
of unbalanced waves under a standard ensemble Kalman filter. 

Speaking on a more abstract level, the MEnK filter also provides a step 
towards a more integrated view on dynamics and data assimilation. The MEnK 
filter is closer to the basic idea of filtering in that regard 
than the current practice of batched assimilation in intervals of 6 to 9 hours, 
which is more in line with four dimensional variational
data assimilation (which is a smoothing technique). The close relation between
our mollified ensemble Kalman filter, continuous time Kalman filtering and nudging should also 
allow for an alternative approach to the assimilation of non-synoptic measurements 
\citep{sr:evensen}.

An IAU ensemble Kalman filter has also been implemented. It is observed that the
filter eventually becomes unstable and/or inaccurate unless a damping factor of
$\gamma = 1.0$ is applied to the fast wave part of the modified Lorenz-96 model. 
We did not attempt to optimize the performance of the IAU assimilation approach. It is feasible 
that optimally chosen weights (see \cite{sr:polavarapu04} for a detailed discussion) could lead to 
a stable and accurate performance similar to  the proposed  MEnK filter. It should be noted that
the MEnK filter requires less model integration steps than the IAU ensemble Kalman filter.
The same statement applies to the recent work of \cite{sr:stauffer09}. Indeed, the
work of \cite{sr:stauffer09} is of particular interest since it demonstrates the benefit of
distributed assimilation in the context of the more realistic shallow-water equations. We would
also like to emphasize that IAU ensemble Kalman filtering, the hybrid approach of \cite{sr:stauffer09}
and our MEnKF filter have in common that they point towards a more ``continuous'' assimilation of
data. Such techniques might in the future complement the current operational practice of batched 
assimilation in intervals of, for example, 6 hours at the Canadian Meteorological Centre 
\cite{sr:houtekamer05}. 

We finally mention that an alternative approach to reducing the generation of unbalanced dynamics 
is to make covariance localization respect balance. This is the approach advocated, for example, by
\cite{sr:kepert09}. \\

\noindent
{\bf Acknowledgment.} We would like to thank Eugenia Kalnay for pointing us
to the work on incremental analysis updates, Jeff Anderson for discussions
on various aspects of ensemble Kalman filter formulations, and Bob Skeel for
clarifying comments on mollification.


\bibliographystyle{plainnat}
\bibliography{survey}


\section*{Appendix}

We have stated in \cite{sr:br10} without explicit proof 
that an ensemble Kalman filter analysis
step at observation time $t_j$ is equivalent to
\begin{equation} \label{proof}
{\bf x}_i(t_j^+) = {\bf x}_i(t_j^-) - \int_0^1 
{\bf P} \nabla_{{\bf x}_i} {\cal V}_j({\bf X})\, {\rm d}s,
\end{equation}
where ${\bf x}_i(t_j^-)$ denotes the forecast values and
${\bf x}_i(t_j^+)$ is the analyzed state of the $i$th ensemble
member. The reader was instead referred to \cite{sr:simon}.

Before we demonstrate how (\ref{proof}) gives rise to (\ref{odeEnKF}),
we provide a derivation of (\ref{proof}) for the simplified
situation of $n=1$, $k=1$, i.e.~the covariance matrices
${\bf P}$ and ${\bf R}$ are scalar quantities and the observation 
operator ${\bf H}$ is equal to one. Under these assumptions the 
standard Kalman analysis step gives rise to
\begin{equation}
P_a = \frac{P_f R }{P_f + R}
\end{equation}
and 
\begin{equation}
\overline{x}_a = \frac{\overline{x}_f R + y P_f}{P_f + R} 
\end{equation} 
for a given observation value $y$. 

We now demonstrate that this update is equivalent to twice the application
of a Kalman analysis step with $R$ replaced by $2R$. Specifically, we obtain
\begin{equation}
P_a' = \frac{2P_m R }{P_m + 2R}, \qquad 
P_m = \frac{2P_f R }{P_f + 2R}
\end{equation}
for the resulting covariance matrix $P_a'$ with intermediate value
$P_m$. The analyzed mean $\overline{x}_a'$ is provided by
\begin{equation}
\overline{x}_a' = \frac{2 \overline{x}_m R + y P_m}{P_m + 2R}, \qquad 
\overline{x}_m = \frac{2 \overline{x}_f R + y P_f}{P_f + 2R} .
\end{equation}
We need to demonstrate that $P_a = P_a'$ and $\overline{x}_a = 
\overline{x}_a'$. We start with the covariance matrix and obtain
\begin{align}
P_a' = \frac{\frac{4P_f R}{P_f + 2R} R}{\frac{2P_f R}{P_f + 2R} + 2R}
= \frac{4P_f R^2}{4P_f R + 4R^2} = \frac{P_f R}{P_f + R} = P_a.
\end{align}
A similar calculation for $\overline{x}_a'$ yields
\begin{align}
\overline{x}_a' = \frac{2 \frac{2\overline{x}_f R + y P_f}{P_f + 2R} R
+ y \frac{2P_f R}{P_f + 2R}}{2R + \frac{2P_f R}{P_f + 2R}}
= \frac{4\overline{x}_f R^2 + 4y P_f R}{4R^2 + 4R P_f} =
\frac{\overline{x}_f R + y P_f}{P_f + R} = \overline{x}_a .
\end{align}
Hence, by induction, we can replace the standard Kalman analysis step by
$K>2$ iterative applications of a Kalman analysis with $R$ replaced by
$K R$. We set $P_0 = P_f$, $\overline{x}_0 = \overline{x}_f$ and perform
\begin{equation}
P_{k+1} = \frac{K P_k R}{ P_k + K R}, \qquad 
\overline{x}_{k+1} = \frac{K \overline{x}_k R + y P_k}{P_k + K R}
\end{equation}
for $k=0,\ldots,K-1$. We finally set $P_a = P_K$ and $\overline{x}_a = 
\overline{x}_K$. Next we introduce a step-size $\Delta s = 1/K$ and assume
$K\gg 1$. Then
\begin{align}
\overline{x}_{k+1} = \frac{\overline{x}_k R + \Delta s y P_k}{R + 
\Delta s P_k} 
= \overline{x}_k - \Delta s P_k R^{-1} \left( \overline{x}_k 
- y \right) + {\cal O}(\Delta s^2)
\end{align}
as well as
\begin{equation}
P_{k+1} = \frac{P_k R}{R + \Delta s P_k} =
P_k - \Delta s P_k R^{-1} P_k + {\cal O}(\Delta s^2).
\end{equation}
Taking the limit $\Delta s \to 0$, we obtain the two differential
equations
\begin{equation}
\frac{{\rm d}P}{{\rm d} s}  = - P R^{-1} P , \qquad
\frac{{\rm d} \overline{x} }{{\rm d} s}  = -P R^{-1} \left( 
\overline{x} - y \right)
\end{equation}
for the covariance and mean, respectively. The equation for $P$ can be 
rewritten in terms of its square root $Y$ (i.e.~$P = Y^2)$ as 
\begin{equation}
\frac{{\rm d} Y}{{\rm d} s}  = -\frac{1}{2} P R^{-1} Y.
\end{equation}
We now consider an ensemble consisting of two members $x_1$ and $x_2$. 
Then $\overline{x} = (x_1 + x_2)/2$ and 
\begin{equation}
P = \frac{1}{2} (x_1-x_2)^2, \qquad
Y = \frac{1}{\sqrt{2}} (x_1-x_2).
\end{equation}
Hence 
\begin{equation}
\frac{\rm d}{{\rm d} s} (x_1 + x_2) = -PR^{-1} ( x_1 + x_2 - 2y),
\qquad
\frac{\rm d}{{\rm d} s} (x_1 - x_2) = -\frac{1}{2} PR^{-1} ( x_1 - x_2),
\end{equation}
which gives rise to
\begin{equation}
\frac{{\rm d} x_i}{{\rm d} s}  = -\frac{1}{2} PR^{-1} \left(x_i + \overline{x} -2y\right)
\end{equation}
for $i=1,2$. Upon integrating these differential equations from $s=0$ with initial conditions
$x_i(0) = x_i(t_j^-)$, $i=1,2$, to $s=1$ and setting $x_i(t_j^+) = x_i(1)$
we have derived a particular instance of the general
continuous Kalman formulation (\ref{proof}). 

The general, non-scalar case
can be derived from a multiple application of Bayes' theorem as follows.
We denote the Gaussian prior probability distribution function by $\pi_f({\bf x})$ and
the analyzed Gaussian probability density function by $\pi_a({\bf x})$. Then
\begin{equation} \label{B1}
\pi_a ({\bf x}) \propto \pi_f({\bf x}) \times \exp (-S_j({\bf x}))
\end{equation}
with $S_j({\bf x})$ defined by (\ref{S}). Bayes' formula (\ref{B1}) can be reformulated
to 
\begin{equation}
\pi_a ({\bf x}) \propto \pi_f({\bf x}) \times \prod_{k=1}^K \exp \left(-\frac{1}{K} S_j({\bf x} \right),
\end{equation}
which gives rise to incremental updates 
\begin{equation}
\pi_{k+1}({\bf x}) \propto \pi_k({\bf x}) \times 
\exp \left(-\frac{1}{K} S_j({\bf x} \right) 
\end{equation}
for $k=0,\ldots,K-1$ with $\pi_0 = \pi_f$. The incremental update in the probability density functions
$\pi_k$ gives rise to a corresponding incremental Kalman update step, which leads to the desired continuous 
formulations in the limit $K\to \infty$.

Having shown that an ensemble Kalman filter analysis step is equivalent
to (\ref{proof}), we next demonstrate that (\ref{odeEnKF}) is the appropriate
mathematical formulation for the complete data assimilation scheme. 
To do so we replace the Dirac delta function $\delta(t-t_j)$ 
in (\ref{odeEnKF}) by a function that takes the constant value
$1/\epsilon$ on the interval $[t_j-\epsilon/2,t_j+\epsilon/2]$ and zero
elsewhere. This function approaches the Dirac delta function as
$\epsilon \to 0$. Then we obtain from (\ref{odeEnKF}) 
\begin{equation}
{\bf x}_i(t_j-\epsilon/2 + \delta) = 
{\bf x}_i(t_j-\epsilon/2) + \int_0^\delta \left[
f({\bf x}_i) - \frac{1}{\epsilon} {\bf P} \nabla_{{\bf x}_i} 
{\cal V}_j({\bf X})\right] {\rm d}\tau 
\end{equation}
for $\delta \in [0,\epsilon]$, 
where we have assumed for simplicity that the ODE (\ref{ode}) is
autonomous. We now rescale integration time
to $s = \tau/\epsilon$, which leads to
\begin{equation}
{\bf x}_i(t_j-\epsilon/2 + \epsilon \Delta) = 
{\bf x}_i(t_j-\epsilon/2) + \int_0^\Delta \left[\epsilon
f({\bf x}_i) -  {\bf P} \nabla_{{\bf x}_i} 
{\cal V}_j({\bf X})\right] {\rm d}s ,
\end{equation}
with $\Delta \in [0,1]$. We arrive at (\ref{proof}) in the limit 
$\epsilon \to 0$ and $\Delta = 1$.

\end{document}